\documentclass[12pt,a4paper]{amsart}

\usepackage{pslatex}
\usepackage{dsfont}
\usepackage{mathrsfs}
\usepackage[colorlinks,pdfdisplaydoctitle]{hyperref}
\usepackage[a4paper,centering]{geometry}
\newtheorem{theorem}{Theorem}[section]

\newtheorem{lemma}[theorem]{Lemma}
\newtheorem{proposition}[theorem]{Proposition}
\theoremstyle{definition}
\newtheorem{definition}[theorem]{Definition}
\theoremstyle{remark}
\newtheorem{remark}[theorem]{Remark}
\numberwithin{equation}{section}
\DeclareMathOperator{\Div}{div}
\DeclareMathOperator{\e}{e}
\newcommand{\R}{\mathbf{R}}
\newcommand{\N}{\mathbf{N}}
\newcommand{\Torus}{\mathds{T}_3}
\newcommand{\Energy}{\mathscr{E}}
\newcommand{\B}{\mathscr{B}}
\newcommand{\F}{\mathscr{F}}
\newcommand{\E}{\mathbb{E}}
\newcommand{\Pb}{\mathbb{P}}
\newcommand{\reg}{\mathcal{R}}
\newcommand{\sreg}{\mathcal{S}}
\newcommand{\cov}{\mathcal{Q}}
\newcommand{\ms}{\mathscr{C}}
\newcommand{\loc}{{\textrm{\tiny loc}}}
\newcommand{\uno}{\mathds{1}}
\newcommand{\scal}[1]{\langle #1\rangle}
\newcommand{\lab}[1]{\textsf{\footnotesize\bfseries [#1]}}
\newcommand{\tu}{\widetilde{u}}
\newcommand{\tv}{\widetilde{v}}
\newcommand{\tz}{\widetilde{z}}
\newcommand{\oomega}{{\overline{\omega}}}
\newcommand{\arxiv}[1]{\href{http://arxiv.org/abs/#1}{\texttt{\small arXiv:#1}}}
\hypersetup{pdftitle={An almost sure energy inequality for Markov solutions to the 3D Navier-Stokes equations},
pdfauthor={M. Romito},
pdfsubject={We show the existence of Markov solutions to the Navier-Stokes equations satisfying the almost sure energy inequality},
pdfcreator={M. Romito},
pdfkeywords={stochastic Navier-Stokes equations, energy inequality, Markov solutions}}
\begin{document}
\title[A.\ s.\ energy inequality for Markov solutions to Navier-Stokes]{An almost sure energy inequality for Markov solutions to the 3D Navier-Stokes equations}
\author[M. Romito]{Marco Romito}
\address{Dipartimento di Matematica, Universit\`a di Firenze, Viale Morgagni 67/a, 50134 Firenze, Italia, {\tt\url{http://www.math.unifi.it/users/romito}}.}
\email{romito@math.unifi.it}
\date{}
\begin{abstract}
We prove existence of weak martingale solutions satisfying an almost sure version of the energy inequality
and which constitute a (almost sure) Markov process.
\end{abstract}
\maketitle
\section{Introduction}

Regardless several attempts, the well-posedness of the martingale problem
for the stochastic (as well as the deterministic case) 3D Navier-Stokes
equations remains an open problem (see \cite{DapDeb08}, \cite{Rom08a}
for details on the martingale problem).

A major breakthrough has been the paper of Da Prato and Debussche~\cite{DapDeb03}
(see also \cite{DebOda06}, \cite{Oda07}, \cite{DapDeb08}), where they show
that there are \emph{special} solutions which correspond to a Markov semigroup
which is strong Feller and uniquely ergodic.

A different approach has been introduced in \cite{FlaRom08} (see also \cite{FlaRom06},
\cite{FlaRom07}, \cite{Rom08}, \cite{Rom08a}) where similar results have been
proved with a completely different method. Here we follow this approach
and consider the Navier-Stokes equations
\begin{equation}\label{e:nse}
\begin{cases}
\dot u - \nu\Delta u + (u\cdot\nabla)u + \nabla p = \dot W,\\
\Div u = 0,
\end{cases}
\end{equation}
with periodic boundary conditions on the 3D torus.

The aim of this paper is to complete the work presented in \cite{Rom08},
where it was proved that the semigroup associated to the Markov solutions
introduced in \cite{FlaRom08} converges to a unique invariant measure.
In order to prove that the rate of convergence is exponential, it was
assumed the existence of Markov solutions satisfying an almost sure
energy inequality. Such result is the main theorem (Theorem~\ref{t:main})
of this paper (see also Remark~\ref{r:link}). The method used here to
prove the main result is essentially the same as in \cite{BloFlaRom08},
where the almost sure energy balance was introduced to handle the
space-time white noise forcing the equation.
\section{Definitions and main result}\label{s:main}

We fix some notations we shall use throughout the paper and we refer to
Temam \cite{Tem95} for a detailed account of all the definitions. Let
$\Torus=[0,2\pi]^3$ and let $\mathcal{D}^\infty$ be the space of infinitely
differentiable divergence-free periodic vector fields $\varphi:\R^3\to\R^3$
with mean zero on $\Torus$. Denote by $H$ the closure of $\mathcal{D}^\infty$
in $L^2(\Torus,\R^3)$ and by $V$ the closure in $H^1(\Torus,\R^3)$.
Denote by $A$, with domain $D(A)$, the \emph{Stokes} operator and
define the bi-linear operator $B:V\times V\to V'$ as the projection
onto $H$ of the nonlinearity of equation~\eqref{e:nse}.
Consider finally the abstract form of problem~\ref{e:nse},
\begin{equation}\label{e:absnse}
du + (\nu Au + B(u,u))\,dt = \cov^{\frac12}dW,
\end{equation}
where $W$ is a cylindrical Wiener process on $H$ and $\cov$ is a linear
bounded symmetric positive operator on $H$ with finite trace. Denote by
$(e_k)_{k\in\N}$ a complete orthonormal system of eigenfunctions
of $\cov$, so that $\cov e_k = \sigma_k^2 e_k$.

Next, we define the probabilistic framework where problem~\eqref{e:absnse}
is considered. Let $\Omega=C([0,\infty);D(A)')$, let $\B$ be the Borel $\sigma$-field
on $\Omega$ and let $\xi:\Omega\to D(A)'$ be the canonical process on
$\Omega$ (that is, $\xi_t(\omega)=\omega(t)$). A filtration can be
defined on $\B$ as $\B_t=\sigma(\xi_s:0\leq s\leq t)$. Let $\Omega^t=C([t,\infty);D(A)')$
and denote by $\B^t$ the Borel $\sigma$-field of $\Omega^t$. Define
the \emph{forward shift} $\Phi_t:\Omega\to\Omega^t$ as
$\Phi_t(\omega)(s)=\omega(s-t)$ for $s\geq t$. Given a probability $P$
on $(\Omega,\B)$, we shall denote by $\omega\mapsto P|^\omega_{\B_t}$,
with $P|^\omega_{\B_t}\in\Pr(\Omega^t)$, a \emph{regular conditional probability
distribution} of $P$, given $\B_t$.

For every $\varphi\in\mathcal{D}^\infty$ consider the process $(M_t^\varphi)_{t\geq0}$
on $\Omega$ defined for $t\geq0$ as
\begin{equation}\label{e:mart_eq}
M_t^\varphi
 = \scal{\xi_t-\xi_0,\varphi}_H
   + \nu\int_0^t\scal{\xi_s,A\varphi}_H\,ds
   - \int_0^t\scal{B(\xi_s,\varphi),\xi_s}_H\,ds.
\end{equation}
\begin{definition}[Weak martingale solution]\label{d:wms}
Given $\mu_0\in\Pr(H)$, a probability $P$ on $(\Omega,\B)$ is a \emph{weak
martingale solution} starting at $\mu_0$ to problem~\eqref{e:absnse} if
\begin{itemize}
\item[\lab{w1}] $P[L^2_\loc([0,\infty);H)]=1$,
\item[\lab{w2}] for each $\varphi\in\mathcal{D}^\infty$ the process $M_t^\varphi$
is square integrable and $(M_t^\varphi,\B_t,P)$ is a continuous martingale
with quadratic variation $[M^\varphi]_t=t|\cov^{\frac12}\varphi|^2_H$,
\item[\lab{w3}] the marginal of $P$ at time $t=0$ is $\mu_0$.
\end{itemize}
\end{definition}
Define for every $k\in\N$ the process $\beta_k(t)=\sigma_k^{-1}M_t^{e_k}$. Under
a weak martingale solution, $(\beta_k)_{k\in\N}$ is a sequence of independent
one dimensional Brownian motions. Thus, under any martingale solution, the process
\begin{equation}\label{e:wiener}
W(t) = \sum_{k=0}^\infty\sigma_k\beta_k(t)e_k
\end{equation}
is a $\cov$-Wiener process and
\begin{equation}\label{e:ou}
z(t) = W(t) - \nu\int_0^t A\e^{-\nu A(t-s)}W(s)\,ds
\end{equation}
is the associated Ornstein-Uhlenbeck process starting at $0$. Define the process
$v(t,\cdot) = \xi_t(\cdot) - z(t,\cdot)$. Since $M_t^\varphi=\langle W(t),\varphi\rangle$
for every test function $\varphi$, it follows that
$$
\scal{v(t)-\xi_0,\varphi}_H + \nu\int_0^t\scal{v(s),A\varphi}_H\,ds - \int_0^t\scal{B(\xi_s,\varphi), \xi_s}_H\,ds = 0.
$$
In different words, under a weak martingale solution $P$, $v$ is a weak solution
(i.\ e.\ a distributional solution) of the equation
\begin{equation}\label{e:Veq}
\dot v + \nu A v + B(v+z,v+z) = 0, \qquad P-\text{a.\ s.},
\end{equation}
with initial condition $v(0)=\xi_0$. An energy balance functional can be associated
to $v$, as it has been done in \cite{FlaRom08} for $\xi$,
\begin{equation}\label{e:Venergy}
\Energy_t(v,z) = \frac12|v_t|_H^2 + \nu\int_0^t|v_r|_V^2\,dr - \int_0^t\scal{z_r,B(v_r+z_r,v_r)}\,dr.
\end{equation}
\begin{definition}[Energy martingale solutions]\label{d:ems}
Given $\mu_0\in\Pr(H)$, a probability $P$ on $(\Omega, \B)$ is a \emph{energy
martingale solution} starting at $\mu_0$ of~\ref{e:absnse} if
\begin{itemize}
\item[\lab{e1}] $P$ is a weak martingale solution (see Definition~\ref{d:wms}),
\item[\lab{e2}] $P[v\in L_\loc^\infty([0,\infty);H)\cap L^2_\loc([0,\infty);V)]=1$,
\item[\lab{e3}] there is a set $T_P\subset (0,\infty)$ of null Lebesgue measure such
	that for all $s\not\in T_P$ and all $t\geq s$,
	$$P[\Energy_t(v,z)\leq\Energy_s(v,z)] = 1.$$
\end{itemize}
\end{definition}
Property~\lab{e3} is meaningful only if the map $\Energy$ is finite, at least almost surely,
with respect to a solution $P$, and measurable. By property~\lab{e2}, $v$ is $P$-a.\ s.
weakly continuous in $H$ and so $|v(t)|_H^2$ is defined point-wise in the energy estimate.
Similarly, the other terms are also $P$-a.\ s.\ finite, again by \lab{e2} and the
regularity properties of $z$ under $P$. Measurability is slightly more challenging
and will be examined later in Proposition~\ref{p:measurable}.
\begin{theorem}\label{t:main}
There exists a family $(P_x)_{x\in H}$ of energy martingale solutions
such that the \emph{almost sure Markov property} holds. More precisely,
for every $x\in H$, for almost every $s\geq0$ (including $s=0$), for all
$t\geq s$ and all bounded measurable $\phi:H\to\R$,
\begin{equation}
\E^{P_x}[\phi(\xi_t')|\B_s] = \E^{P_{\xi_s}}[\phi(\xi_{t-s}')].
\end{equation}
The set of times where the Markov property fails to hold at some point $x$
will be called the set of \emph{exceptional times} of $x$.
\end{theorem}
\begin{remark}\label{r:link}
Under appropriate assumptions on the covariance of the noise, Theorem~$5.12$
of~\cite{FlaRom08} continues to hold thanks to property~\lab{e3}. In particular,
the strong Feller property stated in \cite[Theorem 5.11]{FlaRom08} and the
unique ergodicity in~\cite{Rom08}.

As stated in~\cite{Rom08}, the proof detailed here, together with the corresponding
proof of~\cite[theorem 4.1]{FlaRom08} ensure the existence of the \emph{enhanced
martingale solutions} used for the proof of exponential convergence to the
invariant measure proved in~\cite{Rom08}.
\end{remark}
\section{Proof of Theorem~\ref{t:main}}\label{s:proof}

Before proving the main theorem, it is preliminarily necessary to analyse
more carefully the energy functional~\eqref{e:Venergy}. Let $x\in H$,
and let $z_x(t) = z(t) + \e^{-\nu At}x$ (that is, the solution to the
Stokes problem starting at $x$). Similarly, set $v_x = \xi - z_x$.
\begin{lemma}\label{l:allx}
Let $P$ be an energy martingale solution. Then for every $x\in H$,
$$
P[\Energy_t(v_x,z_x)\leq\Energy_s(v_x,z_x)]=1,
$$
for a.\ e.\ $s\geq0$ (including $s=0$) and every $t\geq s$.
\end{lemma}
\begin{proof}
We just give a sketch of the proof (see for example~\cite{Rom01} for a detailed
proof in a more complicated case). Fix $x\in H$ and set $w(t)=\e^{-\nu At}x$.
Then $z_x=z+w$ and $v_x=v-w$, hence $|v_x(t)|_H^2=|v(t)|_H^2 + |w(t)|_H^2 - 2\scal{v(t),w(t)}_H$,
and, since by assumptions the energy inequality holds for $v$, it is sufficient to
prove a balance equality for $w$ and $\scal{v(t),w(t)}_H$. The balance equality for
$w$ is straightforward by the PDE theory, so we only need to show that for all
$s\geq0$ and $t\geq s$,
$$
\scal{v_t,w_t}_H + 2\nu\int_s^t\scal{v_r,w_r}_V\,dr+\int_s^t\scal{v_r+z_r,B(v_r+z_r,w_r)}\,dr
=\scal{v_s,z_s}_H.
$$
The above formula can be proved by standard methods, since by the regularity of $w$,
$\scal{v_t,w_t}_H$ is differentiable in time and its derivative is $\scal{\dot v,w}+\scal{v,\dot w}$
By replacing $\dot v$ with the corresponding terms in \eqref{e:Veq} and using the
antisymmetric property of the nonlinearity, we get exactly
the above formula.
\end{proof}
\begin{proposition}\label{p:measurable}
Given $x\in H$, the map $(t,\omega)\mapsto\Energy_t(v_x(\omega),z_x(\omega))$,
with $(t,\omega)\in[0,\infty)\times\Omega$, is progressively measurable and
\begin{enumerate}
\item for all $0\leq s\leq t$, the sets
  $E_{s,t}(x)=\{\Energy_t(v_x,z_x)\leq\Energy_s(v_x,z_x)\}$
  are $\B_t$-measurable;
\item for all $t>0$, the sets
  $$E_t(x)=\{\Energy_t(v_x,z_x)\leq\Energy_s(v_x,z_x)\text{ for a.\ e.\ $s\leq t$ (including $0$)}\}$$
  are $\B_t$-measurable;
\item the set
  $$
  E(x)=\reg\cap\{\Energy_t(v_x,z_x)\leq\Energy_s(v_x,z_x)\text{ for a.\ e.\ }s\geq0\text{ (incl. $0$), all }t\geq s\}
  $$
  is $\B$-measurable, where
  \begin{equation}\label{e:regset}
  \reg=\{z\in L^8_\loc([0,\infty);L^4(\Torus)),\ v\in L^\infty_\loc([0,\infty);H)\cap L^2_\loc([0,\infty);V)\}.
  \end{equation}
\end{enumerate}
Moreover, given $P$ satisfying \lab{e1} and \lab{e2}, property \lab{e3} is equivalent
to each of the following:
\begin{itemize}
\item[\lab{e3a}] There is $x\in H$ such that for each $t>0$ there is a set $T\subset(0,t]$ of null
	Lebesgue measure and $P[E_{s,t}(x)]=1$ for all $s\not\in T$.
\item[\lab{e3b}] There is $x\in H$ such that for each $t>0$, $P[E_t(x)]=1$.
\item[\lab{e3c}] There is $x\in H$ such that $P[E(x)]=1$.
\end{itemize}
\end{proposition}
\begin{proof}
Measurability of the map $\Energy$ follows from the semi-continuity properties
of the various term of $\Energy$ with respect to the topology of $\Omega$.
The measurability of each $E_{s,t}(x)$ now follows easily from measurability
of $\Energy$. As it regards sets $E_t(x)$, fix $t>0$ and notice that
the Borel $\sigma$-algebra of the interval $(0,t)$ is countably generated, so
that if $\mathcal{T}_t$ is a countable basis,
$$
E_t(x) = E_{0,t}(x)\cap\bigcap_{T\in\mathcal{T}_t}\{\int_0^t\uno_T(s)(\Energy_t(v_x,z_x)-\Energy_s(v_x,z_x))\,ds\leq0\}
$$
and all sets
$\{\int_0^t\uno_T(s)(\Energy_t(v_x,z_x)-\Energy_s(v_x,z_x))\,ds\leq0\}$
are $\B_t$-measurable by the measurability of $\Energy$.
We next show measurability of $E(x)$. Let $J\subset[0,\infty)$ be a countable dense subset
and define
$$
\reg_t=\{z\in L^8([0,t);L^4(\Torus)),\quad v\in L^\infty(0,t;H)\cap L^2(0,t;V)\},
$$
(notice that the regularity of $z$ and $v$ implies that of $v_x$ and $z_x$),
then $\reg_t\in\B_t$ and, by the lower semi-continuity of the
various terms of $\Energy_t(v_x,z_x)-\Energy_s(v_x,z_x)$
with respect to $t$, it follows that
$$
E(x)=\bigcap_{t\in J}(\reg_t\cap E_t(x))
$$
is $\B$-measurable.
The last statement of the lemma is now obvious from the above equalities, property
\lab{e2} and regularity of $z$.
\end{proof}
We can proceed with the proof of Theorem~\ref{t:main}.
\begin{proof}[\textsc{Proof of Theorem~\ref{t:main}} ]
We use Theorem 2.8 of~\cite{FlaRom08}. For every $x\in H$ let $\ms(x)$
be the set of all energy martingale solutions starting at $\delta_x$.
It is sufficient to show that $(\ms(x))_{x\in H}$ is an \emph{a.\ s.\ pre-Markov
family} (see \cite[Definition 2.5]{FlaRom08}), namely,
\begin{enumerate}
\item each $\ms(x)$ is non-empty, compact and convex, and $x\mapsto\ms(x)$ is
	measurable (with respect to the Borel $\sigma$-field induced by the
	Hausdorff distance on compact sets).
\item For each $x\in H$ and each $P\in\ms(x)$, $P[C([0,\infty);H_\sigma)]=1$, where
	$H_\sigma$ is the space $H$ endowed with the weak topology.
\item For each $x\in H$ and $P\in\ms(x)$ there is a set $T\subset(0,\infty)$ of null
	Lebesgue measure such that for all $t\not\in T$ the following properties hold:
	\begin{enumerate}
	\item (\emph{disintegration}) there is $N\in\B_t$ with $P(N)=0$ such that $\omega\in H$
		and $P|^\omega_{\B_t}\in\Phi_t\ms(\omega(t))$ for all $\omega\not\in N$,
	\item (\emph{reconstruction}) $P\otimes_t Q_\cdot\in\ms(x)$ for each $\B_t$-measurable
		map $\omega\mapsto Q_\omega$ with values in $\Pr(\Omega^t)$ such that there is
		$N\in\B_t$ with $P(N)=0$ and $\omega(t)\in H$, $Q_\omega\in\Phi_t\ms(\omega(t))$
		for all $\omega\not\in N$,
	\end{enumerate}
\end{enumerate}
where $P\otimes_t Q_\cdot$ is the gluing of $P$ and $\omega\mapsto Q_\omega$ (see
\cite[Chapter 6]{StrVar79} for details). 
The first two properties are proved in Lemma~\ref{l:compact}, \emph{disintegration} is proved
in Lemma~\ref{l:disintegrate} and \emph{reconstruction} is proved in Lemma~\ref{l:reconstruct}. 
\end{proof}
\begin{lemma}\label{l:compact}
For each $x\in H$ the set $\ms(x)$ is non-empty compact convex and for all $P\in\ms(x)$,
$P[C([0,\infty);H_\sigma)]=1$. Moreover, the map $x\mapsto\ms(x)$ is Borel measurable.
\end{lemma}
\begin{proof}
For the existence of weak martingale solutions see~\cite{FlaRom08}. The proof of
the energy inequality \lab{e3} is indeed easier than the corresponding martingale
property and can be carried on as in the deterministic case (see for example~\cite{Tem95}
or \cite{Rom01}).
It is easy to show that $\mathcal{C}(x)$ is convex, since
all requirements of Definition~\ref{d:ems} are linear with respect to measures
$P\in\mathcal{C}(x)$. Finally, since under any $P\in\ms(x)$ $z$ is continuous
with values in $H$, weak continuity of $\xi$ follows from property \lab{e2}.

In order to prove both compactness and measurability it is sufficient to show
that for each $x\in H$, for each sequence $(x_n)_{n\in\N}$ converging to $x$
in $H$ and for each $(P_n)_{n\in\N}$, with $P_n\in\ms(x_n)$, the sequence
$(P_n)_{n\in\N}$ has a limit point $P\in\ms(x)$ (this follows from \cite[Lemma 12.1.8]{StrVar79}). 

Let $x_n\to x$ in $H$ and let $P_n\in\ms(x_n)$. The energy inequality for
each $P_n$ ensures that $(P_n)_{n\in\N}$ is tight in $\Omega\cap L^2_\loc([0,\infty);H)$,
hence, up to a subsequence, $P_n\rightharpoonup P_\infty$ for some probability $P_\infty$.
It remains to show that $P_\infty\in\ms(x)$. Property \lab{e1} can be proved essentially
as in \cite[Lemma 4.3]{FlaRom08} (without the complicacies of the super-martingale
property).

Before proving \lab{e2}, \lab{e3}, we use Skorokhod theorem: there are a
probability space $(\Sigma, \F, \Pb)$ and random variables $(\tu_n,\tv_n,\tz_n)_{n\in\N}$,
$(\tu_\infty,\tv_\infty,\tz_\infty)$ such that each $(\tu_n,\tv_n,\tz_n)$ has the
same law of $(\xi,v,z)$ under $P_n$ for $0\leq n\leq\infty$ and $\tz_n\to\tz_\infty$
in $L^8_\loc([0,\infty);L^4(\Torus))$, $\tv_n\to\tv_\infty$ in $\Omega\cap L^2_\loc([0,\infty);H)$
and weakly in $L^2_\loc([0,\infty);V)$.

Property \lab{e2} now follows by semicontinuity of the norms of spaces $L^\infty(0,T;H)$
and $L^2(0,T;V)$ with respect to the topology where $\tv_n\to\tv_\infty$, for all $T>0$.
In view of Proposition~\ref{p:measurable}, we finally prove~\lab{e3a} (with $x=0$).
Fix $T>0$, then since, by the proof of \lab{e1}, $\E\|\tv_n-\tv_\infty\|^2_{L^2(0,T;H)}\to0$,
there is a null Lebesgue set $S\subset (0,T]$ such that for all $s\not\in S$,
$$
\Pb[|\tv_{n'}(s)|_H\to|\tv_\infty(s)|_H\text{ for a subsequence }(\tv_{n'})_{n'\in\N}]=1.
$$
Notice that $0\not\in S$, since we already know that $\tv_n(0)\to\tv_\infty(0)$.
We are now able to prove \lab{e3a} (with $x=0$) for $P_\infty$. For each $n\in\N$
there is a null Lebesgue set $T_n\subset(0,T]$ such that 
$\Pb[\Energy_t(\tv_n,\tz_n)\leq\Energy_s(\tv_n,\tz_n)]=1$, for all $s\not\in T_n$.
Let $T_\infty = S\cup\bigcup_n T_n$ and consider $s\not\in T_\infty$. We have that
$\Energy_t(\tv_n,\tz_n)\leq\Energy_s(\tv_n,\tz_n)$ for all $n\in\N$ and, in the limit
as $n\to\infty$, by virtue of the convergence informations on $\tv_n$ and $\tz_n$ and
of the semicontinuity properties of norms,
$$
\Energy_t(\tv_\infty,\tz_\infty) - \Energy_s(\tv_\infty,\tz_\infty)
\leq\liminf_n (\Energy_t(\tv_n,\tz_n) - \Energy_s(\tv_n,\tz_n)),
$$
and in conclusion \lab{e3a} is true.
\end{proof}
Given $s\geq0$ and $x\in H$, denote by $z(t|s,x)$ the Ornstein-Uhlenbeck process
starting in $x$ at  time $s$, namely
$$
z(t|s,x)=\e^{-\nu A(t-s)}x + (W(t)-W(s)) - \nu\int_s^t A\e^{-\nu A(t-r)}(W(r)-W(s))\,dr.
$$
In particular, $z(t|0,x) = z_x(t)$.  Set moreover $v(t|s,x) = \xi - z(t|s,x)$.
Given $t_0>0$, it is easy to see that for all $\omega\in\Omega^{t_0}$,
$
W(t,\Phi_{t_0}^{-1}(\omega)) = W(t+t_0,\omega) - W(t_0,\omega),
$
and it depends only on the values of $\omega$ in $[t_0, t_0+t]$. Similarly,
\begin{equation}\label{e:shifted}
\begin{aligned}
z(t, \Phi_{t_0}^{-1}(\omega)|s, x) = z(t + t_0, \omega|s + t_0, x),\\
v(t, \Phi_{t_0}^{-1}(\omega)|s, x) = v(t + t_0, \omega|s + t_0, x).
\end{aligned}
\end{equation}
\begin{lemma}\label{l:disintegrate}
For every $x\in H$ and $P\in\ms(x)$, there is a set $T\subset(0,\infty)$ of null
Lebesgue measure, and for all $t\not\in T$ there is $N\in\B_t$,
with $P[N]=0$, such that for all $\omega\not\in N$, $\omega(t)\in H$
and $P|^\omega_{\B_t}\in\Phi_t\ms(\omega(t))$.
\end{lemma}
\begin{proof}
Fix $x\in H$ and $P\in\ms(x)$, let $T_P$ be the set of exceptional times of $P$
and fix $t_0\not\in T_P$. We shall look for a $P$-null set $N\in\B_{t_0}$, with
$N=N_1\cup N_2\cup N_3$, such that $\omega(t_0)\in H$ and
$P|^\omega_{\B_{t_0}}\in\Phi_{t_0}\ms(\omega(t_0))$.

The proof of property \lab{e1} is the same as the proof of Lemma $4.4$ of~\cite{FlaRom08}
and it provides a $P$-null set $N_1\in\B_{t_0}$ out of which \lab{e1} holds.

For any interval $J\subset[0,\infty)$, set $S_J = L^\infty_\loc(J;H)\cap L^2_\loc(J;V)$.
In order to prove \lab{e2}, we need to show that
$P^\omega_{\B_{t_0}}[v(\cdot|t_0,0)\in S_{[t_0,\infty)}]=1$. Set
\begin{equation}\label{e:vdisreg}
\begin{aligned}
\sreg_{t_0} &= \{v\in S_{[0,t_0]},\ \e^{-\nu A(t-t_0)}z(t_0)\in S_{[0,\infty)}\},\\
\sreg^{t_0} &= \{v(\cdot|t_0,0)\in S_{[t_0,\infty)}\},
\end{aligned}
\end{equation}
then $\sreg_{t_0}\in\B_{t_0}$, $\sreg^{t_0}\in\B^{t_0}$ and, from property~\lab{e2},
it follows that $\sreg_{t_0}\cap\sreg^{t_0}$ is a $P$-full set, since
$v(t|t_0,0) = v(t) + \e^{-\nu A(t-t_0)}z(t_0)$. By disintegration,
$$
1=P[\sreg_{t_0}\cap\sreg^{t_0}]=\E[\uno_{\sreg_{t_0}}P|^\omega_{\B_{t_0}}[\sreg^{t_0}]],
$$
hence there is a $P$-null set $N_2\in\B_{t_0}$ such that $P|^\omega_{\B_{t_0}}[\sreg^{t_0}]=1$
for all $\omega\not\in N_2$.

Finally, we prove that \lab{e3c} holds for conditional probabilities. Let $\reg_{t_0}$
be defined as $\reg$ in~\eqref{e:regset} but on the time interval $[0,t_0]$, and set
\begin{align*}
A       &= \reg\cap\{\Energy_t(v,z)\leq\Energy_s(v,z)\text{\ for a.\ e.\ }s\geq0\text{\ (including $0$, $t_0$), all }t\geq s\}\\
A_{t_0} &= \reg_{t_0}\cap\{\Energy_t(v,z)\leq\Energy_s(v,z)\text{ for a.\ e.\ }s\in[0,t_0]\text{ (incl. $0$, $t_0$), all }t\in[s,t_0]\},
\end{align*}
We have that $A_{t_0}\in\B_{t_0}$ and, since $t_0\not\in T_P$, $P[A]=P[A_{t_0}]=1$. Now, if
$\oomega\in A_{t_0}$, set $B(\oomega)=A\cap\{\omega:\omega=\oomega\text{ on }[0,t_0]\}$, then
$B(\oomega)$ is equal to
$$
\{\Energy_t(v^\oomega,z^\oomega)\leq\Energy_s(v^\oomega,z^\oomega)\text{ for a.\ e.\ }s\geq t_0\text{ (including }t_0\text{), all }t\geq s\},
$$
since $v(t+t_0,\omega)=v(t+t_0,\omega|t_0,z(t_0,\omega))$ (and similarly for $z$),
and we have set $v^\oomega(\cdot)=v(\cdot|t_0,z(t_0,\oomega))$ and
$z^{\oomega}(\cdot)=z(\cdot|t_0,z(t_0,\oomega))$. Moreover, the map
$\omega\mapsto\uno_{A_{t_0}}(\omega)P|^\omega_{\B_{t_0}}[B(\omega)]$
is $\B_{t_0}$-measurable, since $P|^\omega_{\B_{t_0}}[B(\omega)]=P|^\omega_{\B_{t_0}}[A]$.
By \lab{e3c} (with $x=0$) for $P$, and disintegration, 
$1=P[A]=\E^P[\uno_{A_{t_0}}(\cdot)P|^\cdot_{\B_{t_0}}[B(\cdot)]]$,
and so there is $N_3\in\B_{t_0}$ such that $P|^\omega_{\B_{t_0}}[B(\omega)]=1$ for all
$\omega\not\in N3$, hence \lab{e3c} holds, with $x=z(t_0,\omega)$.
\end{proof}
\begin{lemma}\label{l:reconstruct}
For every $x\in H$ and $P\in\ms(x)$, there is a set $T\subset(0,\infty)$
of null Lebesgue measure such that $P\otimes_t Q_\cdot\in\Phi_t\ms(\omega(t))$
for all $t\not\in T$. The statement must hold for each $\B_t$-measurable map
$\omega\mapsto Q_\omega$ with values in $\Pr(\Omega^t)$ such that there is
$N_Q\in\B_t$ with $P[N_Q]=0$, and $\omega(t)\in H$, $Q_\omega\in\Phi_t\ms(\omega(t))$,
for all $\omega\not\in N_Q$.
\end{lemma}
\begin{proof}
Let $x\in H$, $P\in\ms(x)$ and let $T_P$ be the set of exceptional times of $P$.
Fix $t_0\not\in T_P$ and let $\omega\mapsto Q_\omega$ according to the statement
of the lemma. Everything boils down to show that $P\otimes_{t_0} Q_\cdot\in\ms(x)$,
and, as in the proof of the previous lemma, we refer to~\cite{FlaRom08} (see
Lemma $4.5$) for the proof of~\lab{e1}.

To verify \lab{e2}, consider again sets $\sreg_{t_0}$ and $\sreg^{t_0}$
defined in \eqref{e:vdisreg}. By \lab{e2} for $Q_\cdot$, for each
$\omega\not\in N_Q$, $Q_\omega[\sreg^{t_0}]=1$. Moreover, by \lab{e2}
for $P$, it follows that $P[\sreg_{t_0}]=1$. Finally, since we know that
$v(t+t_0,\omega)=v(t+t_0,\omega|t_0,0)-\e^{-\nu At}z(t_0,\omega)$,
it follows easily that $\sreg_{t_0}\cap\sreg^{t_0}=\{v\in S_{[0,\infty)}\}$
and so
$$
(P\otimes_{t_0}Q_\cdot)[v\in S_{[0,\infty)}]
 = (P\otimes_{t_0}Q_\cdot)[\sreg_{t_0}\cap\reg^{t_0}]
 = \E^P[\uno_{\sreg_{t_0}}Q_\omega[\sreg^{t_0}]]
 = 1.
$$

Finally, we prove \lab{e3c}. Define $A$ and $A_{t_0}$ as in the proof of
the previous lemma. Since $t_0\not\in T_P$ and $A_{t_0}\in\B_{t_0}$, we know
that $(P\otimes_{t_0} Q_\cdot)[A_{t_0}]=P[A_{t_0}]=1$.
Define $B(\oomega)=A\cap\{\omega:\omega=\oomega\text{ on }[0,t_0]\}$
and notice that, if $\omega\in A_{t_0}\cap N_Q^c$ (which is again a
$\B_{t_0}$-measurable $(P\otimes_{t_0} Q_\cdot)$-full set), then by \lab{e3c}
(with $x=z(t_0,\omega)$) for $Q_\omega$ it follows that $Q_\omega[B(\omega)]=1$.
The map $\omega\mapsto\uno_{A_{t_0}\cap N_Q^c}(\omega)Q_\omega[B(\omega)]$
is then trivially $\B_{t_0}$-measurable and equal to $1$, $P$-a.\ s..
Moreover, $Q_\omega[A]=Q_\omega[B(\omega)]=1$ for all $\omega\in A_{t_0}\cap N_Q^c$
and so
$$
(P\otimes_{t_0} Q_\cdot)[A]
 = \E^P\bigl[\uno_{A_{t_0}\cap N_Q^c}Q_\cdot[B(\cdot)]\bigr]
 = P[A_{t_0}\cap N_Q^c]
 = 1,
$$
in conclusion, \lab{e3c} (with $x=0$) holds for $P\otimes_{t_0} Q_\cdot$.
\end{proof}

\end{document}